\documentclass{article}
\usepackage{amssymb,amsmath,amsfonts,latexsym}
\begin{document}
\title{A Riesz-Haviland type result for truncated moment problems with solutions in $L^1$}
\author{C.-G. Ambrozie}
\date{}

\newtheorem{remark}{Remark}
\newtheorem{proposition}[remark]{Proposition}
\newtheorem{examples}[remark]{Examples}
\newtheorem{example}[remark]{Example}
\newtheorem{lemma}[remark]{Lemma} 
\newtheorem{corollary}[remark]{Corollary}
\newtheorem{theorem}[remark]{Theorem}

\maketitle
\begin{abstract}

We give a  version of the Riesz-Haviland theorem for truncated  moments problems, characterizing the existence of the representing measures that are absolutely continuous with respect to the Lebesgue measure.
The existence of such representing densities describes the dense interior of the convex cone of all data having nonnegative Borel representing measures. A natural regularity assumption on the support is required.

Keywords: positive functional, moments problem, measure

MSC-class: 44A60 (Primary) 49J99 (Secondary)
\end{abstract}

\section{Introduction}

In this paper we  consider  problems of moments in $n$ real variables 
$t=(t_1 ,\ldots t_n )$, with respect to a finite number of  monomial functions $t^i \! =\! t_{1}^{i_1}\cdots t_{n}^{i_n}$ where $i\! =\! (i_1 ,\ldots ,i_n )\in (\mathbb{Z}_{+})^n$, for Lebesgue integrable  densities $f\! =\! f(t )\! \geq \! 0$ a.e. (almost everywhere) on  closed subsets $T\subset \mathbb{R}^n$.
Given a set of numbers $g_{i}$ with $|i|=i_1 +\cdots +i_n \leq 2k$ where $k\in \mathbb{N}$ is fixed, the truncated problem of moments, called also  $T$-problem of moments \cite{Curto} when $T$ is prescribed, is concerned  with
 the existence of the Borel measures $\nu \geq 0$ on $\mathbb{R}^n$ supported on $T$ such that $\int_T t^i d\nu (t)=g_i $ for all $i$ \cite{1}, \cite{19}. One calls $\nu$ a representing measure for $g$, and $g_i$ the moments of $\nu$. We are interested in those  representing measures $\nu =f\, dt $, called {\it representing densities} of $g$, that are absolutely continuous with respect to the
Lebesgue measure $dt =dt_1 \cdots dt_n$.  For any subset $I\subset \mathbb{Z}_{+}^n$,
let $P_I$ denote the linear span of the monomials $X^i$ ($i\in I$) in
 $ \mathbb{R}[X_1 ,\ldots ,X_n ]$. 
 Given $g=(g_i )_{i\in I}$, the linear {\it Riesz functional} $\varphi_g :P_I \to \mathbb{R}$ associated to  $g$ is defined as it is known  \cite{cjot} by
 $\varphi_g X^i =g_i$ for $i\in I$. We say that $\varphi_g$ is $T$-positive \cite{RH}   if $\varphi_g p\geq 0$ for all polynomials $p$ such that $p(t)\geq 0$ for all $t\in T$. 
 This condition is necessary for $g$ to have  representing measures $\nu$ on $T$, since in this case $\varphi_g p=\int_T pd\nu $ for all polynomials $p$.
We  remind below the  Riesz-Haviland theorem \cite{RH}, a basic result concerned with the full problem of moments when $I=\mathbb{Z}_{+}^n$.
\begin{theorem}
\label{rh}
{\em \cite{RH}} Let $T\subset \mathbb{R}^n$ be  closed, and $g=(g_i )_{i\in \mathbb{Z}_{+}^n}$ be a set of reals. 
Then $g$ has representing measures on $T$ if and only if $\varphi_g $ is $T$-positive.
\end{theorem}

An analogue of the Riesz-Haviland theorem for the truncated case was established by R.E. Curto and L.A. Fialkow  \cite{Curto}. For  $I=I_{2k} =\{ i:|i|\leq 2k\} $, when $P_{I}$ is the space of all polynomials of degree $\leq 2k$, they proved that a set $g=(g_i )_{|i|\leq 2k}$ has  representing measures on $T$ if and only if the Riesz functional $\varphi_g :P_{I_{2k}}\to \mathbb{R}$ admits
$T$-positive extensions $\tilde{\varphi}_{g}:P_{I_{2k+2}}\to \mathbb{R}$ to the space $P_{I_{2k+2}}$ of all polynomials of degree $\leq 2k+2$.
 
 By  Theorem \ref{poz}, we characterize the existence of the  representing densities of $g$ on regular supports $T$ by the condition: $\varphi_g p>0$ for all  $p\in P_{I}\setminus \{ 0\}$ such that $p(t)\geq 0$ for all $t\in T$. 
A similar characterization   \cite{junk}  holds under more specific hypotheses, in particular 
if a distinguished moment $\tau_{i_0}$ ($=t^{i_0}$ or linear combination of the $t^i$'s) exists such that  $\lim_{\| t\| \to \infty}\frac{|t^i |}{1+\tau_{i_0} (t)} =0$ for $i\not =i_0$. We mention also the existence of similar results    \cite{LFN} for quadratic ($k=1$) and some quartic ($k=2$, $n=2$) $T$-problems of moments. 

By Theorem \ref{interior}, the set of all $g$ having representing densities is the  dense interior of the set of all $g$ having  representing measures. Our proofs are independent of the results from \cite{Curto}, \cite{junk} and rely mainly on Theorem \ref{rh} \cite{RH}.

 Started by works of Stieltjes, Hausdorff,  Hamburger and Riesz, the area of the truncated problems of moments 
 saw interesting development in various other directions, that we do not attempt to cover. A few recent works
\cite{CSS} - \cite{cufi} should be mentioned in this sense, see also \cite{7}, \cite{11},  \cite{KM}, \cite{14''}. 

I express  thanks to professor Raul Curto and professor 
Lawrence Fialkow for drawing the analogue of the Riesz-Haviland 
theorem in the truncated case  \cite{Curto} to my attention.

\section{Main results}

{\bf Definitions} Let $T\! \not \! =\! \emptyset$ be a closed subset of $\mathbb{R}^n$, such that for any $t\in T$ and $\varepsilon >0$ the Lebesgue measure of the set $\{ x\in T: \| x-t\| <\varepsilon\}$ is $>0$. This always holds a.e. (in the density points of $T$ [7.6, \cite{Rudd}]), but we require it in every point $t$. We call such a  $T$ {\it regular}. Here $\| \, \|$ is the usual Euclidian norm. For any multiindex $i\in \mathbb{Z}_{+}^n$, set $\sigma_i = \{ \! j\! \in \!\mathbb{Z}_{+}^n \! :\! j_k \! =\! \mbox{either}\, 0\, \mbox{or}\, i_k ,\, 1\leq k\leq n\}$.
Let  $I\subset \mathbb{Z}_{+}^n$ be  finite, $I\! \not \! =\!\emptyset$ such that $\sigma_i \subset I$ for all $i\in I$. 
We call such an  $I$ {\em regular}, too. 
Let $g=(g_i )_{i\in I}$ be a set of real numbers with $g_0 =1$. 
We call a convex cone $C$ in a real linear space $F$  {\it acute} if $C\cap (-C)=\{ 0\}$. 
Given  also a linear functional $f:F\to \mathbb{R}$, we  write  $f >0$ if $fc >0$ for all $c\in C\setminus \{ 0\}$.

Lemma \ref{hbs} and its Corollary \ref{extensii} follow from various well known arguments 
(Hahn-Banach, Krein, Mazur, Choquet \cite{Choq}, \cite{DuSc}, \cite{Rudd}) 
on the extension of positive functionals. We found it easier to gather them in a short proof. 
\begin{lemma}\label{hbs}
Let $C\subset F$ be an acute closed convex cone in a finite dimensional  linear space $F$. Let $L\subset F $ be  a linear subspace of codimension 1, 
 and $\phi $ a  linear functional on $L$ such that $\phi \, l>0$ for all $l\in L\cap C$ with $l\not =0$. 
Then there is an extension  $\Phi $ of $\phi$ to $F$ such that $\Phi x>0$ for all $x\in C$, $x\not =0$.
\end{lemma}

{\it Proof}. The sum $Y+\mathbb{R}v$ of a closed convex cone $Y\subset F$ and %
a 1-dimensional subspace $\mathbb{R} v$  is closed. Indeed, if both $\pm v\in Y$, $\mathbb{R}v\subset Y$ and so $Y+\mathbb{R}v =Y$;
if not, we may  suppose  $v\not \in Y$ by replacing $v$ by $-v$ if necessary.
Let $x=\lim_{k\to \infty }(y_k +\lambda_k v)$ where $y_k \in Y$. If  $(\lambda_k )_k$ is bounded,  by compactness we obtain a number $\lambda $  and vector $y\in Y$ such that $x=y+\lambda v$. If it is not, we can assume  either $\lim_{k}\lambda_k =\infty$, or  $\lim_{k}\lambda_k =-\infty$. From $y_k +\lambda_k v \to x$ we derive $ \frac{1}{\lambda_k }y_k +v\to 0$. The case $\lambda_k \to -\infty$ is impossible since it leads to $v=\lim_k \frac{1}{-\lambda_k}y_k \in Y$.  If  $\lambda_k \to \infty$, we obtain $-v=\lim_k \frac{1}{\lambda_k}y_k \in Y$, and so the distance $d(x,Y)$ from $x$ to $Y$ satisfies $d(x,Y)=\lim_k d(y_k -\lambda_k (-v),Y)=\lim_k d(y_k ,Y)=0$ whence $x\in Y$. 
In particular we get (inductively) that the convex cone $K:=C+\ker \phi $ is  closed.

To find  $\Phi$, 
we may suppose $C\not = \{ 0\}$. If $\phi =0$,  $L\cap C =\{ 0\}$. 
Let then $f$ be a linear functional on $F$ with $\ker f=L$. Since $C$ is  acute,  $C\setminus \{ 0\}$ is segmentwise connected. Hence the set $I:=f(C\setminus \{ 0\} )$ is  connected, and so, an interval, that cannot contain $0$ for: $fc=0, c\in C \Rightarrow  c\in L\cap C =\{ 0\} $. Then either $I\subset (0,\infty )$, in which case we let $\Phi =f$, or $I\subset (-\infty ,0)$ in which case $\Phi :=-f$. Then $\Phi c \geq 0$ for all $c\in C$, with strict inequality  if $ c\not =0$. 

If $\phi \not \! =\! 0$, let $F'\! =\! F/\ker \phi$ and $\lambda :F\!\to\! F'$ be the  factorization map. Then $C':=\lambda (C)$,  $L':=\lambda (L)$  and the map $\phi '$ induced by $\phi$ on $L'$ satisfy the  hypotheses as well. 
Indeed, $K\! =\! C\! +\! \ker \lambda$ is closed, $\lambda$ open and $\lambda (F\setminus K)\! =\! F'\setminus C'$ whence $C'$ is closed. Also $\phi \! >\! 0, C\, $acute $ \not\! =\! \{ 0\}$  $\Rightarrow$ $\phi '\! >\! 0,C'\, $acute $\not\! =\! \{ 0\}$. Since $\dim F'\! =\! 2$, we easily get  
the existence of an extension $\Phi '\! >\!0$ of $\phi '$ to $F'$, that will provide $\Phi \! :=\! \Phi '\circ \lambda$: note that
$\dim L'=1$ and $\phi '|_{L'}$ is injective, and so increasing along a  direction of  $L'$ given by a vector $\overline{e}\in \mathbb{R}^2$ (a drawing would help). Then  $ L'=\{ r \overline{e}\}_{r\in \mathbb{R}}$ and $\phi '\overline{e}>\phi ' 0=0$. 
The (closed) convex cone ${\mathcal C}:=\mbox{co} (C',\overline{e})$ generated by $C'$ and $\overline{e}$ 
is acute, for otherwise $\phi ' |_{L'}$ would  decrease and  be $<0$ along the half-line $\{ r\overline{e}\}_{r<0}$ opposite to $\overline{e}$ and contained into $C'$, which  is impossible since $\phi '>0$ on $L'\cap C' \setminus \{ 0\}$.    Since ${\mathcal C}$ is acute, there is  an extension  $\Phi '\! >\! 0$ of $\phi '$, whose kernel is a supporting line for  ${\mathcal C}$  in $0$ only:
fix $\overline{f}\in \mathbb{R}^2$ such that $(\mathbb{R}\overline{f})\cap {\mathcal C}=\{ 0\}$, then
for any $v=r\overline{e}+s\overline{f}$ with $r,s$ real, set $\Phi 'v=r\phi '\overline{e} $. If $v\in 
 C' \setminus \{ 0\}$,   $r> 0$ and so $\Phi 'v=r\phi '\overline{e}>0$.
$\Box$

\begin{corollary}  
  \label{extensii}
  Let $F$ be a finite dimensional  linear space and $K\subset F$ an acute closed convex cone.
  Let  $f_0 :L \to \mathbb{R}$ be a linear functional on a linear subspace $L$ of $F$ such that $f_0 x>0$ for every  $x\in L \cap K$ with $x\not =0$.
  Then there is a linear extension $f$ of $f_0$ to  $F$ such that $f x>0$ for every  $x\in K$, $x\not =0$.
  \end{corollary}
 
For any $\sigma$-finite measure $\mu \geq 0$ on $T$ and $1\leq p\leq \infty$ the symbols  $L^p (T,\mu )$, $L_{+}^p (T, \mu )$ have the usual meaning.
Lemma \ref{feas} is an extension of [Theorem 2.9, \cite{BLl}]. We give its original proof  adapted to our slightly different  context.
  
\begin{lemma}
\label{feas}
{\em (see \cite{BLl})} Let $T\subset \mathbb{R}^n$ be  closed  with positive Lebesgue measure, finite or not. Let $\rho :T\to (0,\infty )$ be locally integrable, and $\mu =\rho dt$ be the measure on $T$ with density $\rho$. Let $I\subset \mathbb{Z}_{+}^n $ be finite. Let $f\geq 0$ a.e. on $T$ be measurable,  $f\not \equiv 0$ a.e.,   such that $\int_T |t^i |f(t)d\mu (t)<\infty$ $(i\in I)$. There is an $r>0$ such that for any $\beta
=(\beta_i )_{i\in I} $ with $\|  \beta \| <r$, there exists a $g=g_\beta \in  L^\infty (T,\mu )$, $g>0$ a.e. with the properties $\int_T |t^i |g(t)d\mu (t)<\infty$,
$$
\int_T t^i g(t)d\mu (t)=\int_T t^i f(t)d\mu (t) +\beta_i
\mbox{ }\mbox{ }\mbox{ }(i\in I).$$
\end{lemma}

{\it Proof}.  Set $T_l \! =\! \{ t\in T\! :\!  f(t)\! \geq\!  1/l,\, \| t\| \leq l\}$.
Using  $\{ t\in T\! :\!  f(t)\! >\! 0\}\! =\! \cup_{l\in \mathbb{N}}T_l $  we find  $\delta >0$ and $T_* \subset T$ bounded with  $0<\mu (T_* )<\infty$ such that $f(t)\geq \delta$ a.e.  on $T_*$. 
The  map $A:L^\infty (T_* )\to \mathbb{R}^N$ ($N=\mbox{card}\, I$) given by  $Au=( \int_{T_*}u(t) t^i d\mu (t))_{i\in I} $
is surjective, for  if there is a vector $\lambda =(\lambda_i )_i \not =0$ orthogonal to its range,  $\sum_i \lambda_i \int_{T_*}u t^i  d\mu  =0$ $\forall \, i,\, u$
whence $\sum_i \lambda_i t^i =0$ a.e. on $T_*$, that is impossible because
  the set of zeroes of a polynomial $ \not \! \equiv\! 0$ has measure zero.
   Then $A$ is open. Hence $0$ is in the interior $\mbox{int}\, C$ of the set $C\! =\! \{ Au\! :\!  u\! \in\!  L^\infty (T_* ,dt), \| u\| \! <\! \delta /4\}$. Fix  $r>0$ such that the ball of center $0$ and radius $r$ is contained in $C$.
   Define  $f_k $  on $T$ by $f_k (t)\! =\! \mbox{min}\, (f(t),k) \! +\! \frac{1}{k}e^{-\| t\| }/\rho (t)$.
   Then $0\leq |t^i |f_k \leq |t^i |f+|t^i |e^{-\| t\|}/\rho (t)\in L^1 (T,\mu )$ and $f_k \to f$ a.e. as $k\to \infty$.
 The vector $(\int_T (f_k   -f ) t^i d\mu  )_{i\in I} \to 0$ in $\mathbb{R}^N$    as $k\to \infty$, by Lebesgue's theorem of dominated convergence.
 Since $0\in \mbox{int}\, C$, for   large $k$  we have $(\int_T (f_k -f) t^i d\mu  )_i \in C$.
   Then $\int_T (f_k -f) t^i d\mu   =\int_T u t^i d\mu  $  for some $u\in L^\infty (T)$ with $\| u\| <\delta /4$ and $u=0$  outside $T_*$. Also, if $\| \beta \| <r$, $\beta \in C$ and so there is a $v\in L^\infty (T )$, $v=0$ outside $T_*$, with $\| v\| <\delta /4$ such that $\int_{T}v t^i d\mu  =\beta_i $. 
   Set $g=f_k -u+v$ for a sufficiently large fixed $k$ ($\geq \delta$). Hence $\int_T g t^i d\mu   =\int_T (f_k -u+v)  t^i d\mu  =\int_T (f_k -u)  t^i d\mu  +\int_T v t^i d\mu  =\int_T f t^i d\mu  +\beta_i$. 
   Since $\rho$ is locally integrable and $T_*$ contained in a ball, on which $u$, $v$ and all $t^i$  are bounded, the functions $t^i u$, $t^i v\,
   \in L^1 (T,\mu )$. Hence $|t^i |g =|t^i |(f_k -u+v)$ is in $L^1 (T,\mu )$. Moreover $f_k , u, v$ (and hence, $g$) are in $L^\infty (T)$.
   On $T_*$, $f\geq \delta $  and $|u|,\, |v|\leq \delta /2$, whence $g=f_k -u+v\geq \mbox{min}\, (f,k) -u+v \geq \delta -u+v \geq \delta /2$. 
    Outside $T_*$, $g=f_k \geq \frac{1}{k}e^{-\| t\|}/\rho (t)>0$. Then $g>0$ a.e. $\Box$
\vspace{3 mm}

\noindent {\bf Notation} Given any closed subset  $T\subset \mathbb{R}^n$    and  finite subset $I\subset \mathbb{Z}_{+}^n$,
  set
  
\noindent   $
  \Gamma_{TI} \! =\! \{ \gamma\! =\! (\gamma_i )_{i\in I}:\,  \exists \mbox{ } \mbox{\rm  Borel measures } \nu \geq 0\, \mbox {\rm  on}\, T\, \mbox{\rm with}\,
  \int_T t^i d\nu (t)\! =\! \gamma_i ,\,  i\in I\} 
  $
  
  and
 
 \noindent  $
  G_{TI}\! =\! \{ g\! =\! (g_i )_{i\in I}\not \! =\! 0:  \, \exists \mbox{ }\, f\in L_{+}^1 (T,dt)\, \mbox{ }\mbox{\rm such that}\,
  \int_T t^i f(t)dt\! =\! g_i ,\, i\in I\} ,
  $
  
  where $|t^i | $ are implicitely supposed to be integrable. 
   
\begin{lemma}\label{densitate}
 Let $T\! \subset\! \mathbb{R}^n$ be a closed regular set and $I\! \subset\!  \mathbb{Z}_{+}^n$  a finite regular set. 
 Then $G_{TI}$ is dense in $\Gamma_{TI}$.
\end{lemma}  

{\it Proof}.
Let $\gamma \in \Gamma_{TI}$. There is a measure $\nu \! \geq \! 0$ on $T$ such that 
 \begin{equation}\label{unuuu}\int_{T} y^i d\nu (y)=\gamma_i \mbox{ }\mbox{ }\mbox{ (}i\in I\mbox{)},\end{equation}
 in particular $\nu (T)<\infty $ since $0\in I$. 
For any   $\varepsilon \in (0,1)$ let $h_{\varepsilon}$ be the characteristic function of the ball $b_\varepsilon$  of center $0$ and radius $\varepsilon$ in $\mathbb{R}^n$, $h_\varepsilon = 1$ on $b_\varepsilon $ and $h_\varepsilon = 0$ outside $b_\varepsilon$. 
For any $y\in T$, let $v_\varepsilon (y)$ be the $n$-dimensional volume of the set $\{ x\in T: \| x-y\| <\varepsilon \}$. Then $v_\varepsilon =(h_{\varepsilon}*h_{T})|_T$ is the convolution of  $h_{\varepsilon}$ with the characteristic function $h_T$ of $T$. Hence the map $y\mapsto v_{\varepsilon}(y)$ is measurable. All $v_\varepsilon (y)>0$  since $T$ is regular.
For $t\in T$, set $\nu_\varepsilon (t)=\int_T \frac{1}{v_\varepsilon (y)}h_\varepsilon (t-y) d\nu (y).$
  By the Tonelli and Fubini theorems, $\nu_\varepsilon\in L_{+}^1 (T,dt)$  has finite moments of orders $i\in I$ on $T$, that we compute by
  \begin{equation}\label{aprox}
\int_T t^i \nu_\varepsilon (t)dt=\int_T  \frac{1}{v_{\varepsilon}(y)}\int_T t^i h_\varepsilon (t-y)dt \, d\nu (y)=\int_T \psi_{\varepsilon i} (y) d\nu (y)   
\end{equation}
 where $\psi_{\varepsilon i} (y)=\frac{1}{v_{\varepsilon}(y)}\int_T t^i h_\varepsilon (t-y)dt$. By the change of variables $t-y=w$, \begin{equation}\label{conv}\psi_{\varepsilon i}  (y)=\frac{1}{v_\varepsilon (y)}\int_{\| w\| <\varepsilon ,w\in T-y}(y+w)^i dw =y^i \frac{1}{v_\varepsilon (y)}\int_{\| w\| <\varepsilon ,w\in T-y} dw  +E_i (\varepsilon, y)$$
 where $E_i (\varepsilon ,y)$ is a linear combination with  binomial  coefficients $c_{ij}$, $$E_i (\varepsilon, y)=\sum_{0\leq j\leq i, j\not =i}c_{ij}y^j \frac{1}{v_\varepsilon (y)}\int_{\| w\| <\varepsilon ,w\in T-y}w^{i-j}dw\end{equation}
 and the order on $\mathbb{Z}_{+}^n$ is given as usual by $j\leq i\, \Leftrightarrow \, j_k \leq i_k$ for $1\leq k\leq n$.
 The set $\{\| w\| <\varepsilon ,w\in T-y\}$ is taken by the translation $w\mapsto w+y$ into the set $\{ x\in T: \| x-y\| <\varepsilon \}$ the Lebesgue measure of which is $v_\varepsilon (y)$.
Hence \begin{equation}\label{dirac}
\frac{1}{v_\varepsilon (y)}\int_{\| w\| <\varepsilon ,w\in T-y} dw 
=1.\end{equation} Then \begin{equation}\label{concl}\psi_{\varepsilon i} (y) =y^i +E_i (\varepsilon , y).\end{equation}
Since  $\| w\| <\varepsilon <1$ and $|i-j|\geq 1$ for all $j$ in (\ref{conv}), $| w^{i-j}| \leq \| w\|^{|i-j|}\leq \| w\| <\varepsilon$. Hence  by (\ref{conv}) and (\ref{dirac}), we obtain the estimate  \begin{equation}\label{estimare}|E_i (\varepsilon ,y)|\leq \varepsilon \sum_{0\leq j\leq i} c_{ij}|y^j |.\end{equation} Since $I$ is regular,  $j\leq i\, \Rightarrow \, |y_{k}^{j_k}|\leq |y_{k}^{i_k}|+1\Rightarrow  |y^j |\leq \prod_{k=1}^n (|y_{k}^{i_k}|+1)=\sum_{\iota \in \sigma_i }|y^\iota |\leq \sum_{\iota \in I }|y^\iota |$ and so we can integrate in (\ref{concl}), (\ref{estimare}) with respect to $\nu$ on $T$. By (\ref{aprox}) and (\ref{unuuu}), this gives $\lim_{\varepsilon \to 0}\int_{T} t^i \nu_\varepsilon (t)dt = \gamma_i$ for all $i\in I$. Set $\tilde{\gamma }_\varepsilon =(\int_{T} t^i (\nu_\varepsilon (t)+\varepsilon e^{-\| t\|})dt )_{i\in I}$. 
Then $\tilde{\gamma}_\varepsilon \in G_{TI}$ and $\lim_{\varepsilon \to 0}\tilde{\gamma}_\varepsilon =\gamma $. 
$\Box$

  \begin{theorem}
  \label{interior}
  Let $T\!\subset\! \mathbb{R}^n$ be a closed regular set, and $I\!\subset \!\mathbb{Z}_{+}^n$  a finite regular set. 
  Then $G_{TI}$ is the dense interior of $\Gamma_{TI}$.
  \end{theorem}
  
  {\it Proof}. 
  By Lemma \ref{feas} for $\rho \equiv 1$,  $G_{TI}$ is open, and so  contained in the interior of $\Gamma_{TI}$ 
  (the regularity of $T$, $I$ is not required here). 
 Let $\gamma$ be in the interior of $\Gamma_{TI}$. 
 There is an $r>0$ such that the ball $B$ of center $\gamma$ and radius $r$ is contained in $\Gamma_{TI}$ 
 (a drawing will be helpful).
By Lemma \ref{densitate}, there is a $\tilde{\gamma} \in G_{TI}\cap B$. 
  By Lemma \ref{feas}  applied to a representing density $f$ of $\tilde{\gamma}$,
  there is an $r'_0 >0$ such that all  balls $B(\tilde{\gamma},r ')$ of  center $\tilde{\gamma}$ and radii $r ' \in (0,r'_0 ]$ satisfy $B(\tilde{\gamma},r ')\subset G_{TI}$.
  We can fix an $r'$ sufficiently small so  that $B(\tilde{\gamma},r ')\subset B$. Let $\gamma '$ be the unique point  such that $\gamma =\frac{1}{2}(\tilde{\gamma}+\gamma ')$. Since $B(\tilde{\gamma},r ')\subset B$,  
  then $B(\gamma ',r')\subset B$; a quick argument to this aim is that $B$ is symmetric with respect to its center $\gamma$ and   $B(\gamma ',r')$, $B(\tilde{\gamma},r ')$ are symmetric to each other over $\gamma$.
In particular $\gamma '\in B\subset \Gamma_{TI}$. By Lemma \ref{densitate}, there is a $\tilde{\gamma}' \in G_{TI}\cap B(\gamma ' ,r')$. Since  $\tilde{\gamma}' \in B(\gamma ' ,r')$, the  point $v$ such that $\gamma =\frac{1}{2}(\tilde{\gamma}' +v)$ must be in $B(\tilde{\gamma},r ')$ ($\subset G_{TI}$). Hence $v\in G_{TI}$.
  Since $G_{TI}$ is  convex, $\gamma =\frac{1}{2}(\tilde{\gamma}' +v)$ and both $\tilde{\gamma}', \, v\in G_{TI}$,
 then $\gamma \in G_{TI}$. 
  $\Box$
\vspace{3 mm}
    
By  the previous results, the following completion (Theorem \ref{poz}) can be made to the Riesz-Haviland theorem \cite{RH} and its truncated version \cite{Curto}.

Condition (b) from below is equivalent to the existence of  a  constant $c>0$ such that $\varphi_g p \geq c\| p\|$ for every $p\in P_I$ with $p(t)\geq 0$ for all $t\in T$, where $\| \, \|$ is any norm on $P_I$: use the compactness of  $\{ p\in P_I :p\geq_{\stackrel{\mbox{ }}{T}} 0, \| p\| =1\} $ and write $\varphi_g (p/\| p\| )>0$ for such $p\in P_I \setminus \{ 0\}$. 
 \begin{theorem}
 \label{poz}
 Let $T\subset \mathbb{R}^n$ be a closed regular set. Let $I\subset \mathbb{Z}_{+}^n$ be a finite regular set. Let $g=(g_i )_{i\in I}$ be a set of numbers with $g_0 =1$. The following statements are equivalent:
 
{\em (a)} There exist functions  $f\in L_{+}^1 (T,dt)$ such that $\int_T |t^i |f(t)dt<\infty$ and
 $$
 \int_T t^i f(t)dt=g_i \mbox{ } \mbox{ for all }i\in I
 ;$$
 
{\em (b)} The Riesz functional $\varphi_g$ defined on the linear span $P_I \subset  \mathbb{R}[X_1 ,\ldots ,X_n ]$ of the monomials $X^i$, $i\in I$ by $$\varphi_g \sum_{i\in I} c_i X^i =\sum_{i\in I} g_i c_i$$
 satisfies  $\varphi_g \, p>0$ for every  $p\in P_I \setminus \{ 0\}$ such that $p(t)\geq 0$ for all $t\in T$.
 \end{theorem} 
  
  {\em Proof}. (a) $\Rightarrow$ (b) 
From $\varphi_g p=\int_T pfdt$ we obtain  as usual that  $\varphi_g p \geq 0$ for every polynomial $p$ such that $p(t)\geq 0$ for all $t\in T$. If moreover $p\not =0$ in $\mathbb{R}[X_1 ,\ldots ,X_n ]$,
the set $Z=\{ t\in \mathbb{R}^n :p(t)=0\}$ of the  zeroes of $p$ is an algebraic variety, (empty or) of dimension $\leq n-1$, and so has null Lebesgue measure.
Then $\varphi_g p>0$, for the equality $\int_T pf\, dt =\int_{T\setminus Z} pf\, dt=0$ with $p>0$ on $T\setminus Z$ would compel $f=0_{L^1 (T,dt)}$ 
that is impossible since $\int_T fdt =g_0 =1$.
  
  (b) $\Rightarrow$ (a) Endow $P_I $   with a norm
  and its dual $P_{I}^*$ with the dual norm. Let $C$ denote the convex cone of all $p\in \mathbb{R}[X_1 ,\ldots ,X_n ]$ such that $p(t)\! \geq \! 0$ for all $t\in T$. Since the Lebesgue measure of $T$ is $>0$, $C$ is acute.
  There is a  constant $c\! =\! c_g \! >\! 0$ such that $\varphi_g p \! \geq\! c\| p\|$ for 
every $p\in P_I \cap C$, see the comment just before the Theorem. 
  Since the map $ \gamma  \mapsto \varphi_\gamma \in P_{I}^*$ is linear,  there is a constant $c'\! >\! 0$ such that $\|\varphi_\gamma \|\! \leq\! c' \| \gamma \|$ for all $\gamma$.   Then for every $\tilde{g}$ in the ball $B$ of center $g$ and         radius $r\! =\! c/\, 2c'$, we have
   $\varphi_{\, \tilde{g}\, } p \! \geq\!  (c/2)\| p\|$ for all $p\in P_I \cap C$. Indeed, 
  $\varphi_{\tilde{g}} p =\varphi_g p +\varphi_{\tilde{g}-g}p\geq c\| p\| -|\varphi_{\tilde{g}-g}p|\geq c\| p\| -c'\| \tilde{g}-g\|\, \| p\| \geq (c/2)\| p\|$.
   Hence condition (b), briefly $\varphi_{\tilde{g}}>0$,  holds as well for all $\tilde{g}$ in the neighborhood $B$ of $g$.
  Write $\mathbb{Z}_{+}^n =\cup_{l\in \mathbb{N}}I_l$ as an increasing union of finite subsets $I_l $  such that  $I_1 =I$.
Let  $P_l$ be the linear span of the monomials $X^i$  with $i\in I_l$. Thus $P_1 \! =\! P_I$. Let $\tilde{g}\! =\! (\tilde{g}_i )_{i\in I} $ be an arbitrary point in $B$, that for the moment we fix. Since $\varphi_{\tilde{g}}>0$, by applying successively Corollary \ref{extensii} for $F\! =\! P_{l+1}$, $K\! =\! C\cap P_{l+1}$ and $L\! =\! P_l$ with $l\! \geq \! 1 $ we obtain, starting from $\psi_{1\tilde{g}} :=\varphi_{\tilde{g}} :P_1 \! \to\!  \mathbb{R}$, a sequence of linear functionals $\psi_{l\tilde{g}} : P_{l}\to \mathbb{R}$ such that  $\psi_{l+1\, \tilde{g}} |_{P_{l}}\! =\! \psi_{l\tilde{g}}$ for all $l\! \geq\! 1$ and $\psi_{l\tilde{g}} p\! >\! 0$ for all $p\in P_l \cap C \setminus \{ 0\}$.
Then we have a  linear functional  $\psi_{\tilde{g}} :\mathbb{R}[X_1 ,\ldots ,X_n ]\! \to\!  \mathbb{R}$, extending  $\varphi_{\tilde{g}}$, determined by  $\psi_{\tilde{g}} |_{P_l}\! =\! \psi_{l\tilde{g}}$ for all $l\geq 1$. 
For any polynomial $p\in C$, $p\not \! \equiv \! 0$  
   there is an $l\! \geq \! 1$
 such that $p\in P_l$, and so $\psi_{\tilde{g}} p\! =\! \psi_{l\tilde{g}} p>0$. Hence $\psi_{\tilde{g}}  p\geq 0$ for every polynomial $p$ such that $p(t)\geq 0$  for all $t\in T$. By the Riesz-Haviland's Theorem \ref{rh}  there is a  measure $\nu_{\tilde{g}}\geq 0$ on $T$, with finite moments of any order, such that $\psi_{\tilde{g}} p\! =\! \int_T pd\nu_{\tilde{g}}$ for every polynomial $p$. In particular for $p\! :=\! X^i$ with $i\in I$, we obtain 
 $\tilde{g}_i \! =\! \varphi_{\tilde{g}} X^i \! =\! \psi_{\tilde{g}} X^i \! =\! \int_T t^i d\nu_{\tilde{g}} (t)$. Thus $ \tilde{g}_i$ ($i\in I$) are the moments of a measure $\nu_{\tilde{g}}$ on $T$, that is, $\tilde{g}\in \Gamma_{TI}$. Since $\tilde{g}$ was arbitrary in a neighborhood of $g$, it follows that $g$ is in the interior of $\Gamma_{TI}$. Then by Theorem \ref{interior}, $g\in G_{TI}$. 
  $\Box$
\vspace{3 mm}

{\bf Acknowledgement} The work was supported by the grants IAA100190903  GAAV and 201/09/0473 GACR, RVO: 67985840.

Institute of Mathematics, AS CR

Zitna 25

115 67 Prague 1

Czech Republic
\vspace{1 mm}

{\it ambrozie@math.cas.cz}
\vspace{5 mm}

and: Institute of Mathematics "Simion Stoilow" - Romanian Academy, 

\hspace{8 mm} PO Box 1-764, 014700 Bucharest, Romania

  \end{document}